\documentclass[12pt]{amsart}

\usepackage{url}

\usepackage{amssymb}
\usepackage[all]{xy}
\usepackage[T1]{fontenc}
\usepackage[utf8]{inputenc}
\usepackage{amsmath, amssymb, amsthm}
\usepackage{tikz}

\DeclareMathOperator{\Id}{Id}
\newcommand{\id}[1]{\Id(#1)}

\newtheorem{thm}{Theorem}[]

\newtheorem{cor}[thm]{Corollary}
\newtheorem{ex}[thm]{Example}

\begin{document}

%
%
%


\author[Berele-Cerulli Irelli-De Loera Chavez-Pascucci]{Allan Berele\and Giovanni Cerulli Irelli \and Javier De Loera Chávez\and Elena Pascucci}
\address{Department of Mathematical Sciences, De Paul University, Chicago, IL 60614}
\email{aberele@depaul.edu}
\address{Dipartimento SBAI, Sapienza Universit\`a di Roma, Via Scarpa 10, Roma (IT) 00161}
\email{giovanni.cerulliirelli@uniroma1.it}  
\email{javieralejandro.deloerachavez@uniroma1.it} 
\email{elena.pascucci@uniroma1.it}

\title{Polynomial identities for quivers via incidence algebras}

\begin{abstract} 
We show that the path algebra of a quiver satisfies the same polynomial identities of an algebra of matrices, if any. In particular, the algebra of $n\times n$ matrices is PI-equivalent to the path algebra of the oriented cycle with $n$ vertices. 
\end{abstract}
\maketitle

Let $Q$ be a finite quiver with vertex set $Q_0$ and let $F$ be a field of characteristic zero. Recall that a path of length $\ell$ in $Q$ is a concatenation of $\ell$ composable arrows of $Q$. Given a path $p$ in $Q$ we denote by $s(p)$ its starting vertex and by $t(p)$ its terminal vertex. For every vertex $i$ of $Q$ we denote by $e_i$ the zero length path with $s(e_i)=t(e_i)=i$. The path algebra of $Q$ is the $F$-algebra linearly spanned by all paths of $Q$, and it is denoted by $FQ$. The product of two paths $p$ and $q$ in $Q$ is zero if $t(p)\neq s(q)$ and it is the path $pq$ otherwise. We notice that $FQ$ is a unitary associative algebra with $1=\sum_{i\in Q_0}e_i$.

Let $\pi$ be an arbitrary set of paths in $Q$. We denote by $FQ_\pi$ the subalgebra of $FQ$ generated by $\pi$. Thus, $FQ_\pi$ is the linear span of the set $\tilde{\pi}$ of all paths in $Q$ that are (non-empty) products of elements of $\pi$. 

In this paper we prove that $FQ_\pi$  has the same polynomial identities of an incidence algebra in the sense of \cite{B}.

For a given positive integer $n$, we denote by $M_n(F)$ the $F$-algebra of $n\times n$ matrices with entries in $F$. We denote by $e_{ij}$ the elementary $n\times n$ matrix whose entries are all zero except the entry $(i,j)$ which is equal to one.  

Given a quiver $Q$ with vertex set $Q_0=\{1,2,\cdots, n\}$, we define the following linear map on the path basis of $Q$
\[
\varphi_Q: FQ\rightarrow M_n(F): p\mapsto e_{s(p)t(p)}.
\]
We denote by $A_Q:=\mathrm{im}(\varphi_Q)$ its image. Thus 
\[A_Q=\{A\in M_n(F)\mid A_{i,j}=0\textrm{ if there are no paths in }Q\textrm{ from }i\textrm{ to }j\}.\]

More generally, we define  $A_\pi:=\varphi_Q(FQ_\pi)$. Thus,
\[A_\pi=\{A\in M_n(F)\mid A_{i,j}=0\textrm{ if there are no paths in }\tilde{\pi}\textrm{ from }i\textrm{ to }j\}.\]
The algebra $A_\pi$ is the incidence algebra, as defined in \cite{B}, associated to the transitive order on the set $Q_0$ given by  $i<j$ if there is a path in $\tilde{\pi}$ from $i$ to $j$. Note that every incidence algebra in the sense of \cite{B} is of the form $A_\pi$.

For an $F$-algebra $A$ we denote by $\id{A}$ the $T$-ideal\footnote{A $T$-ideal is a two-sided ideal of  $F\langle X\rangle$ invariant under all the endomorphisms of $F\langle X\rangle$. Here $F\langle X\rangle$ denotes the free object in the category of algebras containing $A$; thus it is unital if $A$ is. 
}  of polynomial identities satisfied by $A$. The algebra $A$ is called $PI$ if $\id{A}\neq \{0_{F\langle X\rangle}\}$. In \cite{CDP} it is proved that $FQ$ is a PI algebra if and only if $Q$ has no vertices adjacent to more than one oriented cycle and such quivers are called PI. (Note that if there is a vertex $i$ of $Q$ and two distinct paths $p$ and $q$ in $Q$ with both endpoints equal to $i$, then they generate a subalgebra of $FQ$ isomorphic to the free algebra in two generators, and hence $FQ$ is not PI.) On the other hand, polynomial identities of incidence algebras are described in \cite{B}. In this paper we relate these two classes of algebras. Our main result is the following.

\begin{thm}\label{Thm:Main}
If $Q$ is $PI$ then $\id{FQ_\pi}=\id{A_\pi}$ for any set $\pi$.
\end{thm}

If $\pi$ is the set of all paths in $Q$ then $A_\pi=A_Q$ and hence we get:
\begin{cor}
If $Q$ is $PI$ then $\id{FQ}=\id{A_Q}$.
\end{cor}
There are  many examples of infinite dimensional algebras $A$ such that $\id{A}=\id{M_n(F)}$; for instance, one can take $A=M_n(R)$, where $R$ is a commutative infinite-dimensional $F$-algebra. 
As a corollary of Theorem~\ref{Thm:Main} we obtain another example of this kind, but of a different nature. 
Note that for $Q=C_n$, the oriented cycle with $n$ vertices, $A_Q=M_n(F)$. Thus we obtain

\begin{cor}
For any $n\ge 2$, $\id{FC_n}=\id{M_n(F)}$.
\end{cor}

\section{Proof of Theorem~\ref{Thm:Main}}

Since $A_\pi$ is a quotient of $FQ_\pi$, $\id{FQ_\pi}\subseteq \id{A_\pi}$. (This holds for any quiver $Q$, not necessarily PI).

Let us prove that $\id{A_\pi}\subseteq \id{FQ_\pi}$. Let $f(x_1,\cdots, x_m)\in \id{A_\pi}$. Since $F$ has characteristic zero, we can assume that $f$ is multilinear and it has the form 
\[
f(x_1,\cdots, x_m)=\sum_{\sigma\in S_m} \lambda_\sigma x_{\sigma(1)}x_{\sigma(2)}\cdots x_{\sigma(m)},
\] 
where $S_m$ denotes the symmetric group on $m$ elements.
Since $FQ_\pi$ is linearly generated by paths in $\tilde{\pi}$ and $f$ is multilinear, it is enough to show  that $f(p_1,\cdots, p_m)=0$ for every choice $\beta=(p_1,\cdots, p_m)$ of paths in $\tilde{\pi}$. Without loss of generality, we can assume that $p_1\cdots p_m\neq 0$ and hence $p_1\cdots p_m$ is a path in $Q$. 
We have
\[
f(p_1,\cdots p_m)=\sum_{i,j\in Q_0}e_i\left(\sum_{\sigma\in (S_\beta)_{i,j}} \lambda_\sigma p_{\sigma(1)}p_{\sigma(2)}\cdots p_{\sigma(m)}\right)e_j
\]
where $(S_\beta)_{i,j}$ is the set of all permutations $\sigma\in S_m$ with the property that $s(p_{\sigma(1)}p_{\sigma(2)}\cdots p_{\sigma(m)})=i$
and $t(p_{\sigma(1)}p_{\sigma(2)}\cdots p_{\sigma(m)})=j$. Since $Q$ is PI, by \cite[Lemma~6.3]{CDP}, $p_{\sigma(1)}p_{\sigma(2)}\cdots p_{\sigma(m)}=p_{\sigma'(1)}p_{\sigma'(2)}\cdots p_{\sigma'(m)}$ for every $\sigma,\sigma'\in (S_\beta)_{i,j}$. With the shorthand $p_{i,j}^\beta:=p_{\sigma(1)}p_{\sigma(2)}\cdots p_{\sigma(m)}$, we have
\[
f(p_1,\cdots p_m)=\sum_{i,j\in Q_0}\left(\sum_{\sigma\in (S_\beta)_{i,j}} \lambda_\sigma\right) p_{i,j}^\beta.
\]
Then, since $\varphi_Q$ is a homomorphism of algebras and $f\in \id{A_\pi}$, 
\[
0=f(\varphi_Q(p_1),\cdots, \varphi_Q(p_m))=\varphi_Q(f(p_1,\cdots p_m))=\sum_{i,j\in Q_0}(\sum_{\sigma\in (S_\beta)_{i,j}} \lambda_\sigma ) e_{i,j}
\]
and hence 
\[
\sum_{\sigma\in (S_\beta)_{i,j}} \lambda_\sigma=0
\]
for all $i,j=1,\cdots n$. This concludes the proof.

\section{Examples}
Let us illustrate Theorem~\ref{Thm:Main} in a few examples. Following \cite{B}, we denote by $T_0$ a singleton set with no relations, and by $T_n$ the set with $n$ elements with the relation $T_n\times T_n$. For every set $S$ endowed with a transitive relation $<$,  we denote by $I(S)$ the T-ideal of the incidence algebra of $S$. Note that $I(T_n)=\id{M_n(F)}$.

\begin{ex}
Let $Q=\xymatrix{1\ar^\alpha[r]&2}$ be a quiver of type $A_2$ and let $\pi=\{e_1,\alpha\}$. Then $A_\pi=\left(\begin{smallmatrix}\ast&\ast\\0&0\end{smallmatrix}\right)$. The incidence algebra $A_\pi$ is associated to the transitive order on $S=\{1,2\}$ given by $1<1$, $1<2$ and thus a chain. The chain decomposes as $S=T_1T_0$ and hence, by Theorem~\ref{Thm:Main} and \cite[Corollary~5]{B},  
\[\id{FQ_\pi}=\id{A_\pi}=I(T_1)I(T_0)=\langle [x_1,x_2]\rangle_T\langle x_3\rangle_T=\langle [x_1,x_2]x_3\rangle_T\] (the last equality follows by Jacobi's identity).
\end{ex}
\begin{ex}\label{Ex1}
Let $Q=\begin{tikzpicture}[>=latex, baseline]
\node (1) at (0,0) {$1$};
\node (2) at (1,0) {$2$};
\node (3) at (2,0) {$3$};
\draw[->] (1) to[out=160, in=200, loop] node[left] {$\gamma$} (1);
\draw[->] (1) to node[above] {$\alpha_1$} (2);
\draw[->] (2) to[out=120, in=60, loop] node[above] {$\beta$} (2);
\draw[->] (3) to node[above] {$\alpha_2$} (2);
\draw[->] (3) to[out=20, in=-20, loop] node[right] {$\delta$} (3);
\end{tikzpicture}$. Then $A_Q=\left(\begin{smallmatrix}\ast&\ast&0\\0&\ast&0\\0&\ast&\ast\end{smallmatrix}\right)$ is the incidence algebra associated with the transitive order on $S=\{1,2,3\}$ given by $1<1<2<2$, $3<3<2<2$ which is a union of two maximal chains. By Theorem~\ref{Thm:Main} and \cite[Theorem~7]{B}, $\id{FQ}=\id{A_Q}=I(\{1,2\})\cap I(\{2,3\})$. Since $I(\{1,2\})=I(\{2,3\})=\id{\left(\begin{smallmatrix}\ast&\ast\\0&\ast\end{smallmatrix}\right)}$, by Malcev's theorem we get $\id{FQ}=\langle [x_1,x_2][x_3,x_4]\rangle_T$.
\end{ex}

\begin{ex}
Let $Q$ be the quiver of Example~\ref{Ex1}
and let $\pi$ be the set of all paths in $Q$ but $\{e_3,\delta\}$. Then $A_\pi=\left(\begin{smallmatrix}\ast&\ast&0\\0&\ast&0\\0&\ast&0\end{smallmatrix}\right)$ is the incidence algebra associated with the transitive order on $S=\{1,2,3\}$ given by $1<1<2<2$, $3<2<2$ which is a union of two maximal chains. By Theorem~\ref{Thm:Main} and \cite[Theorem~7]{B}, $\id{FQ_\pi}=\id{A_\pi}=I(\{1,2\})\cap I(\{2,3\})$. Since $I(\{1,2\})=\id{\left(\begin{smallmatrix}\ast&\ast\\0&\ast\end{smallmatrix}\right)}$, by Malcev's theorem we get $I(\{1,2\}=\langle [x_1,x_2][x_3,x_4]\rangle_T$. Since $\{2,3\}=\{3\}\{2<2\}=T_0T_1$, by \cite{B},  $I(\{2,3\})=I(T_0)I(T_1)$. We get \[\id{FQ_\pi}=\langle [x_1,x_2][x_3,x_4]\rangle_T\cap \langle x_5\rangle_T\langle [x_6,x_7]\rangle_T=\langle [x_1,x_2][x_3,x_4]\rangle_T.\]
\end{ex}

\section{Acknowledgements} GCI, JDLC and EP are supported by NextGenerationEU - PRIN 2022 -B53D23009430006 - 2022S97PMY -PE1- investimento M4.C2.1.1- Structures for Quivers, Algebras and Representations (SQUARE), by “Progetti
di Ateneo” Sapienza Universit\`a di Roma and by GNSAGA-INDAM. We thank the anonymous referee for their useful comments. 

MSC: \textbf{16R10, 16G20}.

\textbf{Keywords}: Polynomial Identities, Quivers, Path Algebras


\begin{thebibliography}{99}\label{bib}
\bibitem{B} A.~Berele, \emph{Incidence algebras, polynomial identities, and an $A\otimes B$ 
counterexample}. Comm. Algebra \textbf{12} (1984), no. 1-2, 139–147. 

\bibitem{CDP} G.~Cerulli~Irelli, J.~De~Loera~Chávez, E.~Pascucci, \emph{Quivers with Polynomial Identities}, \url{https://doi.org/10.48550/arXiv.2508.00662}

\end{thebibliography}
\end{document}